\input amstex
\UseAMSsymbols
\documentstyle{amsppt}
\define\r{\bold u}
\define\rz{{\bold u}_{\zeta}}
\define\U{U^0_{\zeta}}
\define\Uz{U_{\zeta}}
\define\K{\left[K\atop \ell\right]}
\TagsOnRight

\magnification=\magstephalf
\NoBlackBoxes
\hsize=30pc

\vsize=42pc

\hoffset=.5truein

\voffset=.75truein

\document

%\refstyle{A}
\vskip.3truein

\topmatter

\title Quantized Hyperalgebra of Rank 1
\endtitle

\rightheadtext{Quantized Hyperalgebra}

\leftheadtext{W. Chin and L. Krop}

\author W. Chin and L. Krop\footnote"\dag"{The second author gratefully acknowledges support of the
University Research Council of DePaul University} \endauthor

\date January 30, 2004\enddate

\affil DePaul University
\endaffil

\address Department of Mathematics, DePaul University, Chicago,
Illinois 60614
\endaddress

\email wchin\@condor.depaul.edu, lkrop\@condor.depaul.edu
\endemail


\abstract We study the algebra $U_{\zeta}$ obtained via Lusztig's `integral' form [Lu 1, 2] of the generic quantum algebra for the Lie algebra $\frak {g=sl}_2$ modulo the two-sided ideal generated by $K^l-1$. We show that $U_{\zeta}$ is a smash product of the quantum deformation of the restricted universal enveloping algebra $\bold u_{\zeta}$ of $\frak g$ and the ordinary universal enveloping  algebra $U$ of $\frak g$, and we compute the primitive (= prime) ideals of $\Uz$. Next we describe a decomposition of $\bold u_{\zeta}$ into the simple $U$- submodules, which leads to an explicit formula for the center and the indecomposable direct summands of $\Uz$. We conclude with a description of the lattice of cofinite ideals of $\Uz$ in terms of a unique set of lattice generators.

\endabstract

\endtopmatter

\head \S 0 Introduction\endhead

G. Lusztig constructed in [Lu 1,2,3] quantum algebras associated to the defining relations of the finite-dimensional semi-simple Lie algebra $\frak g$. The method used there is similar to the one employed by Kostant [Ko] in his construction of the hyperalgebra for $\frak g$. We refer to Lusztig's algebra as the quantum hyperalgebra of $\frak g$.

Let $\frak {g=sl}_2$ be the rank 1 simple Lie algebra. Fix a field $\Bbbk$ of characteristic zero containing a primitive $\ell$- th root of unity $\zeta$ of an odd order. We let $U_q$ stand for the usual generic quantum algebra of $\frak{sl}_2$. We let $\hat \Uz$ denote the quantum algebra associated with $\frak{sl}_2$ as in [Lu 1]. We define $\Uz$ as the quotient of $\hat U_{\zeta}$ modulo the ideal generated by $K^{\ell}-1$. We denote by $\rz$ the Frobenius - Lusztig kernel in $\Uz$ and let $U$ stand for the ordinary enveloping algebra of $\frak{sl}_2$ over $\Bbbk$.

The goal of this paper is to obtain an explicit description of the primitive ideals, the center, blocks and the lattice of cofinite ideals of $\Uz$.

The paper is organized as follows. In Section 1 we show that $\Uz$ is the smash product of the Frobenius - Lusztig kernel and $U$. This feature of $\Uz$ is one that distinguishes this case from the higher rank cases, and informs the structure of $\Uz$. The theory of simple modules and their annihilators is taken up in Section\nolinebreak 2. Here we point out that the prime ideals are primitive. The main result is the existence of the Steinberg - Lusztig factorization [Lu 1] for every (not necessarily finite- dimensional) simple $\Uz$- module and an explicit formula for a primitive ideal. As a technical preliminary we give the presentation of the ``diagonal" part $\U$ of $\Uz$ by generators and relations. The results of this section admit of a generalization to all semi-simple finite-dimensional $\frak g$. This has been carried out in [CK 2]. However, in the present case the proof proceeds along different lines which yield a stronger result.

In Section 3 we describe the decomposition of $\rz$ into a direct sum of simple $U$- modules. The exposition relies heavily on the structure of principal indecomposable modules for $\rz$. In Section 4 we compute the center of $\Uz$. As a consequence of that we show that every indecomposable direct summand (block) of $\Uz$ is of the form $\epsilon\Uz$ where $\epsilon$ is a block idempotent of $\rz$.

Section 5 contains a description of the lattice of cofinite ideals. An important feature of the lattice is its distributivity. Consequently it has a unique set of lattice generators, namely, the meet- irreducible ideals. We show that each such ideal is the annihilator of a simple, Weyl or co-Weyl module, or else the annihilator of injective hull of a simple finite - dimensional module, viewed as comodule for the finitary dual of $\Uz$. For a more precise statement on the structure of these ideals, see Theorem 7.

It should be noted that the complete lattice of two-sided ideals of $U_q$ is known due to [Ba]. His classification asserts in part that a cofinite ideal is a unique product of maximal ideals. This can be seen directly as follows. Let $I$ be such an ideal. Then $U_q/I$ is finite dimensional, hence a semisimple $U_q$- module. Therefore $I=\cap \frak{m}_r$, where the $\frak{m}_r$ run over all maximal ideals containing $I$. But, as $\text{Ext}_U^1(X,Y)=0$ for all finite dimensional $U_q$- modules $X$ and $Y$, a result in [Mo 1] yields $\cap\frak{m}_r=\prod\frak{m}_r$. 

The main theorems of \S 5 are parallel to the just mentioned result of [Ba]. A description of the entire lattice of two-sided ideals remains an open question, which will be addressed in a later paper.
\vskip.2in

We fix some more notation. The comultiplication in $\Uz$ and $U$ will be denoted by $\Delta_{\zeta}$ and $\Delta$, respectively. An element of the form $\displaystyle\left[{K;c}\atop t\right]$ [Lu 1] defined in $U_q$ will be marked by a subscript $q$. Unsubscribed $\displaystyle\left[{K;c}\atop t\right]$ is $\displaystyle\left[{K;c}\atop t\right]\otimes 1$. An unmodified $\otimes$ means $\otimes_{\Bbbk}$. We recall the definition of the quantum Casimir element. This is $c_{\zeta}=\displaystyle{FE+\frac{\zeta K+\zeta^{-1}K^{-1}}{(\zeta-\zeta^{-1})^2}}$. It is straightforward to check that $c_{\zeta}$ lies in the center of $\Uz$ 

\head \S 1 Smash product Representation of $U_{\zeta}$\endhead

In what follows we put $\frak g$ equal to ${\frak sl}_2$ and denote by $\{e,h,f\}$ its standard basis.

By [Lu 3] there exists a unique Hopf algebra map $\text{Fr}:U_{\zeta}\rightarrow U$, called the Frobenius map, specified on generators of $U_{\zeta}$ as follows: $\text{Fr}(E)=\text{Fr}(F)=0,\,\text{Fr}(K)$\linebreak$=1,\,\text{Fr}(E^{\ell})=e$ and $\text{Fr}(F^{\ell})=f$

Let $\r_{\zeta}$ be the subalgebra of $U_{\zeta}$ generated by $E,F$ and $K$. Put $\r^{+}_{\zeta}=\r_{\zeta}\cap \text{Ker}\, \epsilon$. Combining results of [Lu 3] and [An] we know that $\text{Ker}\,\text{Fr}$ is generated by $\rz^{+}$ as a left or right ideal. Moreover, $\text{Fr}$ induces a right $U$- comodule algebra structure on $U_{\zeta}$ via $\rho:U_{\zeta}\rightarrow U_{\zeta}\otimes U$, $\rho=(I\otimes\text{Fr})\Delta_{\zeta}$, where $I$ is the appropriate identity map. We conclude using a couple of Schneider's results [Sch 1,2] that $\r_{\zeta}=U_{\zeta}^{\text{co}U}$. The upshot of these remarks is that we have an exact sequence in the category of Hopf algebras
$$\Bbbk\longrightarrow \r_{\zeta}\longrightarrow U_{\zeta}\overset\text{Fr}\to\longrightarrow U\longrightarrow\Bbbk$$

We want to show that the above sequence splits by an algebra and right $U$- comodule homomorphism $\gamma:U\rightarrow U_{\zeta}$ satisfying $\text{Fr}\circ \gamma=I$

\proclaim{Lemma 1}The mapping $\gamma$ defined on $\frak g$ by
$\gamma(e)=E^{(\ell)},\,\gamma(f)=F^{(\ell)}$ extends to an algebra homomorphism $U\to\Uz$
\endproclaim
\demo{Proof}Let $H=[E^{(\ell)},F^{(\ell)}]$. We must show that $[H,E^{(\ell)}]=2E^{(\ell)}$ and\linebreak 
$[H,F^{(\ell)}]=-2F^{(\ell)}$. To this end we recall that by [Lu 1] we have $H=\left[K\atop\ell \right]+f$, where $f=\dsize\sum_{1\le t\le \ell-1}F^{(\ell-t)}\left[{K;-2(\ell-t)}\atop t\right]E^{(\ell-t)}$. Since $F^{(s)}\left[{K;c}\atop t\right]=\left[{K;c+2s}\atop t\right]F^{(s)}$ it follows $F^{(\ell-t)}\left[{K;-2(\ell-t)}\atop t\right]E^{(\ell-t)}=\left[K\atop t\right]F^{(\ell-t)}E^{(\ell-t)}$. An easy induction on $s$ leads to the formula $$F^{(s)}E^{(s)}=\dsize\frac1{[s]!^2}\prod\limits_{i=1}^s\left(c_{\zeta}-\frac{\zeta^{2(i-1)+1}K+\zeta^{-2(i-1)-1}K^{-1}}{(\zeta-\zeta^{-1})^2}\right)\tag 1-1$$ 
Thus $f$ lies in the subalgebra of $\r_{\zeta}$ generated by $K$ and $c_{\zeta}$ which is centralized by $E^{(\ell)}$ and $F^{(\ell)}$. It follows that 
$$\align[H,E^{(\ell)}]&=\left[\left[K\atop \ell \right],E^{(\ell)}\right]\\&=E^{(\ell)}\left\{\left[{K;2\ell}\atop \ell\right]-\left[K\atop \ell\right]\right\}
\endalign$$
Applying [Lu 2,(g9)] we obtain in $U_{\Cal A}$
$$\left[{K;2\ell}\atop \ell\right]_q=\sum_{0\le j\le \ell}q^{2\ell(\ell-j)}\left[{2\ell}\atop j\right]K^{-j}\left[K\atop{\ell-j}\right]_q$$
which, in view of $\left[{2\ell}\atop j\right]_{\zeta}=0$ for $0<j<\ell$, reduces to 
$$\left[{K;2\ell}\atop \ell\right]=\left[K\atop \ell\right]+\left[{2\ell}\atop \ell\right]=\left[K\atop \ell\right]+2$$
in $U_{\zeta}$. Thus we arrive at $\bigl[H,E^{(\ell)}\bigr]=2E^{(\ell)}$. Similar computations give\linebreak $\bigl[H,F^{(\ell)}\bigr]=-2F^{(\ell)}$.
\qed
\enddemo

\proclaim{Theorem 1}$U_{\zeta}=\r_{\zeta}\# U$ is the smash product of $\r_{\zeta}$ and $U$
\endproclaim
\demo{Proof} As $\gamma$ is an algebra map, it is convolution invertible. In fact, it is also a $U$- comodule map, i.e. $\rho\circ\gamma=(\gamma\otimes I)\circ\Delta$. Since the mappings in the last equation are algebra homomorphisms it suffices to check it on generators which is straightforward. By [DT], $\gamma$ induces $U$- action in $\r_{\zeta}$ via $v\cdot a=\gamma(v_1)a\gamma^{-1}(v_2),\,v\in U,\,a\in \r_{\zeta}$ (where $\Delta v$ is written as $v_1\otimes v_2$). This action gives rise to the smash product $\r_{\zeta}\# U$ which is isomorphic to $U_{\zeta}$ by op.cit.
\qed
\enddemo

\head \S 2 Simple Modules and Primitive Ideals\endhead

The subalgebra $U^0_{\zeta}$ of $U_{\zeta}$ is by definition generated by $K$ and the elements $\left[{K;c}\atop t\right],\,c\in\Bbb Z,\,t\in\Bbb Z_+$ [Lu 1,2]. It is responsible for the weight theory in the category of $U_{\zeta}$- modules. We are interested in finding generators and relations for $U^0_{\zeta}$. The following result appears in [CK 2]. We include its proof for the reader's convenience

\proclaim{Proposition 1} $U^0_{\zeta}$ is commutative algebra generated by $K$ and $\left[K\atop \ell\right]$ subject to the relation $K^{\ell}=1$
\endproclaim
\demo{Proof}We apply the Doi-Takeuchi theory of cleft extensions to $U^0_{\zeta}$. We remark that $U^0_{\zeta}$ is a subcoalgebra of $U_{\zeta}$, and that $\text{Fr}$ sends $U^0_{\zeta}$ to $U^0$, the subalgebra of $U$ generated by $h$. It follows that $\rho$ restricted to $U^0_{\zeta}$ induces a right $U^0$- comodule structure on the latter. We define an algebra homomorphism $\phi\colon U^0\to U^0_{\zeta}$ by sending $h$ to $\left[K\atop \ell\right]$. We observe that $\phi$ splits $\text{Fr}$ a $U^0$- comodule map. 

We conclude by [DT] that $U^0_{\zeta}=A\#_{\sigma}U^0$, where $A=\displaystyle(U^0_{\zeta})^{\text{co}U^0}$. Further, since $\phi$ is an algebra map, the cocycle $\sigma$ is trivial. Also $U^0_{\zeta}$ is commutative, therefore the action of $U^0$ is trivial as well. Thus $U^0_{\zeta}=A\otimes U^0$. It remains to identify $A$. Let $\r_{\zeta}^0=\Bbbk[K]$ be the subalgebra of generated by $K$. We have that $A$ is the subalgebra of coinvariants in $U^0_{\zeta}$, hence $A=\r_{\zeta}\cap U^0_{\zeta}$. Using the PBW - theorem for $U_{\zeta}$ one can see readily that $A=\r_{\zeta}^0$, and the proof is complete.
\qed
\enddemo

\remark {\bf Remark \rom{1}} We sketch a direct proof of the proposition. First of all we have from [Lu 2] that $U^0_{\zeta}$ is spanned by the set $\{K^{\delta}\left[K\atop t\right]|\delta=0,1;\,t\ge 0\}$. The relation $K^{\ell}=1$ enables us to reduce that set to the subset $\{\left[K\atop t\right]|t\ge 0\}$. The last step is a formula of independent interest (cf. [Lu 1, 4.3]).
\endremark

\proclaim{Lemma 1'} Let $m$ be a positive integer written as $m=m_0+\ell m_1,\,0\le m_0\le \ell -1$. Then  $\left[K\atop m\right]=\left[K\atop {m_0}\right]\displaystyle\binom\K{m_1}$, where the second factor is the ordinary binomial expression
\qed
\endproclaim

\remark{\bf Remark \rom{2}} The lemma fails without assumption $K^{\ell}=1$.
\endremark
The Proposition follows readily. For, the sets $\{\left[K\atop a\right]|0\le a\le \ell-1\}$ and\linebreak  $\{\binom\K b|b\ge 0\}$ form bases for $\r_{\zeta}^0$ and the subalgebra generated by $\K$, respectively.

\remark{\bf Remark \rom{3}} According to the Cartier- Kostant- Milnor- Moore theorem (see [Mo, 5.6.5]) $U^0_{\zeta}=\Bbbk G\otimes U(P)$, where $G$ is the set of group-like and $P$ is the Lie algebra of primitive elements in $U^0_{\zeta}$. The Frobenius map restricted to $U(P)$ induces an isomorphism of $U(P)$ with $U$. Hence $\dim P=1$ and there is a unique primitive $d$ such that $\text{Fr}(d)=h$. We want to give an explicit formula for $d$.
$\r^0_{\zeta}=\Bbbk[K]$ has $\ell$ minimal idempotents $e_i,\,0\le i\le \ell-1$. Explicitly $e_m=\frac1\ell(\displaystyle\sum_{j=0}^{\ell-1}\textstyle\zeta^{-mj}K^j)$. It follows by a direct computation that $\Delta_{\zeta}e_m=\displaystyle\sum_{i+j\equiv m}e_i\otimes e_j$, where the $\equiv$ denotes congruence modulo $\ell$. On the other hand one can verify that $\Delta_{\zeta}(\K)=\K\otimes 1+1\otimes\K+\displaystyle\sum_{i+j\ge \ell}e_i\otimes e_j$. 

We now define $d$ by the formula$$d=\K+\frac1\ell(\displaystyle\sum_{m=1}^{\ell-1}me_m)$$
Using the above formulas it is straightforward to check that $d$ is indeed a primitive element. For an alternate treatment see [CK 2]
\qed
\endremark

We now turn to a computation of characters of $U^0_{\zeta}$. Let $X=\text{Alg}(U^0_{\zeta},\Bbbk)$ be the group of algebra homomorphisms under convolution. To every $\chi\in X$ we associate its weight $\lambda=(r,\alpha)$ according to the equalities $\chi(K)=\zeta^r$ and $\chi(\left[K\atop \ell\right])=\alpha$.

Let $\Lambda\simeq \Bbb {Z}_{\ell}\times \Bbbk$ be the set of all weights. We remark that every integer $m$ can be viewed as a weight. For, writing (uniquely) $m=m_0+m'\ell, 0\le m_0\le \ell-1$, we can identify $m$ with the pair $(m_0,m')$. The next lemma has been proved in [CK 2]. We repeat the proof for the sake of completeness.

\proclaim{Lemma 2}Let $\lambda=(r,\alpha)$ and $\mu=(s,\beta)$ be two weights. Then
$$\lambda+\mu=\cases (r+s, \alpha +\beta), &\text{if $r+s< \ell$}\\
                      (r+s-\ell,\alpha +\beta +1),&\text{else}\endcases$$

\endproclaim 
\demo{Proof}Pick two characters $f$ and $g$ of weight $\lambda$ and $\mu$, respectively. Since $K$ is a group-like $(f\ast g)(K)=\zeta^{r+s}$. Further, recall the formula\linebreak[CP, 11.2] $\displaystyle\Delta(\left[K\atop \ell\right])=\sum_{0\le j\le \ell}\left[K\atop \ell  -j\right]K^{-j}\otimes \left[K\atop j\right]K^{\ell-j}$. In view of $K^{\ell}=1$ and $f(\left[K\atop m\right]=\left[r\atop m\right]_{\zeta}$ for every $m<\ell$ and similarly for $g$ it follows readily that $(f\ast g)(\left[K\atop \ell\right])=\alpha+\sum^{'}+\beta$ where $\sum{'}=\sum_{1\le j\le \ell}\zeta^{-(r+s)j}\left[r\atop\ell-j\right]_{\zeta}\left[s\atop j\right]_{\zeta}$. It remains to notice that as $\left[m\atop\ell\right]_{\zeta}=0$ for all $m<\ell$, the identity [Lu, 2(g 9-10)] shows that $\sum{'}=\left[{r+s}\atop \ell\right]$ and the assertion follows from [Lu, 1(3.2)].
\qed
\enddemo

As in the classical case we define a partial order on $\Lambda$ by saying \linebreak$\lambda\le \mu\Longleftrightarrow \lambda=\mu-2n$ for some $n\in\Bbb Z^+$. Let ${\Bbbk}_{\lambda}$ be the 1-dimensional $U^0_{\zeta}$- module of weight $\lambda=(r,\alpha)$, i.e. ${\Bbbk}_{\lambda}=\Bbbk v_{\lambda}$ with $Kv_{\lambda}=\zeta^r$ and $\left[K\atop \ell\right]v_{\lambda}=\alpha v_{\lambda}$. We put $B^+=U^0_{\zeta}U^+_{\zeta}$ and observe that the $U^0_{\zeta}$ action on ${\Bbbk}_{\lambda}$ can be lifted to the $B^+$- action by setting $E^{(n)}v_{\lambda}=0$ for all $n>0$. We define (as usual) the $U_{\zeta}$- module $V(\lambda)$ by $V(\lambda)=U_{\zeta}\otimes {\Bbbk}_{\lambda}$ with the left regular action of $U_{\zeta}$. One can check easily that $F^{(n)}v_{\lambda}$ has weight $\lambda-2n$. By the standard argument $V(\lambda)$ has a unique maximal submodule and a unique irreducible quotient denoted by $L(\lambda)$.

In preparation for the next statement we fix some notation. For a $U$- module $M$ we denote by $M^{\text Fr}$ the $U_{\zeta}$- module obtained via the pull-back along $\text {Fr}$. For a weight $\lambda=(m_0,m_1)$ with $m_1\in\Bbb Z$ we write $\lambda=m$ where $m=m_0+m_1\ell$. If $m_1=0$ we say that $\lambda$ is {\it $\ell$- restricted}. We abbreviate $L(m_0,m_1)$ to $L(m)$. For an $\alpha\in\Bbbk$ we let ${\overline L}(\alpha)$ denote the highest weight $\alpha$ simple $U$- module. For a weight $\lambda$ we denote by $P(\lambda)$ the primitive ideal $$P(\lambda)=\displaystyle\text{ann}_{U_{\zeta}}L(\lambda)$$
If $\lambda=(r,\alpha)$ we set $$\frak {p}(r)=\displaystyle\text{ann}_{\r_{\zeta}}L(r)\quad {\text and}\quad \overline{P}(\alpha)=\text{ann}_U\overline{L}(\alpha)$$

\proclaim{Theorem 2} (i) Every prime ideal of $\Uz$ is primitive

(ii) Every primitive ideal of $U_{\zeta}$ has the form $P(\lambda)$ for some weight $\lambda$. Furthermore, if $\lambda=(r,\alpha)$ then we have $P(\lambda)=\frak{p}(r)U+\r_{\zeta}\overline{P}(\alpha)$

(iii) Every simple $U_{\zeta}$- module $S$ has the Steinberg - Lusztig factorization
$$S\simeq L(r)\otimes \overline {S}^{\text{Fr}}$$
for a suitable restricted $r$ and a $U$- module $\overline S$. Further $L(r)\otimes \overline{L}(\alpha)^{\text Fr}$ is isomorphic to $L(r,\alpha)$.

\endproclaim

First a Lemma.

\proclaim{Lemma 3} Let $N$ be the nilpotent radical of $\r_{\zeta}$. Then $N$ is stable under $U$ and the $U$- module $\r_{\zeta}/N$ is trivial.
\endproclaim
\demo{Proof} By [Lu 1, 7.1] the restriction of $L(r)$ to $\r_{\zeta}$ is a simple module for every restricted $r$. In fact, these exhaust all simple $\r_{\zeta}$- modules which can be seen directly or following the argument of the classical modular case [Cu]. Thus $N=\displaystyle{\cap_{0\le r\le \ell}}\frak{p}(r)$. Now, $U$- invariance of every prime ideal of $\rz$ follows from [Ch 1] or more immediately from [GW, Prop. 1.1] in view of the fact that every prime of $\rz$ is a minimal prime.

Since $\r_{\zeta}$ is generated by $K,E,F$ it suffices to show that $[E^{(\ell)}, X]$ and $[F^{(\ell)}, X]$ lie in $N$ for every generator $X$. This is clear for $X=K$ as $K$ commutes with $E^{(\ell)}$ and $F^{(\ell)}$. Suppose $X=F$. Then we have $[E^{(\ell)}, F]=\left[{K;-\ell+1}\atop 1\right]E^{\ell -1}$. Also $E^{\ell -1}L(r)=0$ for all $r<\ell$. Thus $[E^{(\ell)},F]$ lies in every $\frak{p}(r)$ for $r\ne \ell$. Let $r=\ell$. Recall  that $L(\ell -1)$ is the span of $\{F^iv_0|0\le i\le \ell -1\}$, where $v_0$ is a primitive vector of weight $\ell -1$. It follows that $E^{\ell -1}L(\ell -1)=\Bbbk v_0$. As $\left[{K;-\ell+1}\atop 1\right]v_0=0$, the proof is complete.
\qed
\enddemo

We proceed to the proof of the theorem. (i) For every $U$- invariant ideal $I$ of $\r_{\zeta}$ it is straightforward to check using the smash product structure that $UI=IU$. Hence $NU$ is a nilpotent ideal of $U_{\zeta}$. In view of the above lemma $\displaystyle{U_{\zeta}/NU}$ is semisimple, and hence $NU$ is the Jacobson radical of $U_{\zeta}$.

Now, let $P$ be a prime ideal of $\Uz$. By the opening remark $P\supset NU$. Passing on to $\Uz/{NU}\simeq \rz/N\otimes U\simeq \otimes M_t(U)$ we may assume that $P$ is a prime ideal of $M_t(U)$ for some $t$. It follows readily that $P=M_t(\frak p)$ for a prime ideal $\frak p$ of $U$. Then by [NG] $\frak p$ is primitive, and by [Po] so is $P$.

(ii) Pick $P$, a primitive ideal of $U_{\zeta}$, and a simple module $S$ with $P=\text {ann}_{\Uz}S$  Let $\frak{p}=P\cap \r_{\zeta}$. Then $\frak{p}\supset N$, hence $\frak{p}$ is a prime ideal of $\rz$. For, every ideal of $\rz$ containing $\frak{p}$ is $U$- invariant. Therefore the inclusion $I\cdot J\subset \frak{p},\;I,J$ ideals of $\r_{\zeta}$ implies $IU\cdot JU\subset P$, hence $I$ or $J$ lies in $\frak{p}$

Choose $r$ such that $\frak{p}=\frak{p}(r)$ and let $\pi:U_{\zeta}\longrightarrow \r_{\zeta}/{\frak{p}}\otimes U$ be the map composed of the natural epimorphism $U_{\zeta}\longrightarrow U_{\zeta}/{\frak{p}}U$ and the isomorphism $U_{\zeta}/{\frak{p}}U\simeq \r_{\zeta}/{\frak{p}}\otimes U$ sending $x\# y+\frak{p}U$ to $(x+\frak{p})\otimes y$. We remark that $\r_{\zeta}/{\frak{p}}$ is just $\text{End}_{\Bbbk}L(r)$ because $L(r)$ is absolutely simple. Indeed, the standard argument [Lu 1] for simplicity of $L(r)$ works for every field containing $\Bbb {Q}(\zeta)$. Denoting $\r_{\zeta}/{\frak{p}}$ by $E$ we see that $S$ is an $E\otimes U$- module. Now, the standard argument with minor modifications as e.g. in [St, Lemma 68]  shows that $S= L(r)\otimes \overline S$ for a suitable $U$- module $\overline S$. 

Let $\overline P$ denote the $\text{ann}_U\overline S$. $E\otimes U/{\overline P}$ acts faithfully in $S$ which implies that $E\otimes \overline P$ is the $\text{ann}_{E\otimes U}S$. We want to compute $T:=\pi^{-1}(E\otimes \overline P)$. To this end we notice the equalities $\r_{\zeta}\cap \frak{p}U=\frak{p}$ and $\overline {P}\cap \frak{p}U=(0)$, both of which come directly from the $\Bbbk$- isomorphism $U_{\zeta}\simeq \r_{\zeta}\otimes U$. But then $\pi(\r_{\zeta}\overline P)=E\otimes \overline P$ which shows that $T=\frak{p}U+\r_{\zeta}\overline P$. On the other hand $\pi(P)=\text{ann}_{E\otimes U}S$ and therefore $P=T$. Further, by Duflo's theorem [Du] $\overline {P}=\overline{P}(\alpha)$ which completes the proof of (ii).

For part (iii) we will show that the identity map $x\otimes y\mapsto x\otimes y: L(r)\otimes \overline S\longrightarrow L(r)\otimes{\overline S}^{\text Fr},x\in L(r),\,y\in\overline S$ is a $U_{\zeta}$- isomorphism. This must be checked on the generators. The latter either belong to $\r_{\zeta}$ or to $\gamma(U)$. We'll do both cases.

Recall that $U_{\zeta}$ acts in $L(r)\otimes\overline{L}(\alpha)^{\text{Fr}}$ via the pullback along $\rho=(I\otimes\text{Fr})\Delta_{\zeta}$ with $U_{\zeta}\otimes U$ acting along the factors. Hence for an $a\in\r_{\zeta}$, in view of $\text{Fr}(a)=\epsilon(a)$, we have $\rho(a)(x\otimes y)=ax\otimes y$ as needed. On the other hand for a $\gamma(v),\,v\in U$ we have $\rho(\gamma(v))=\gamma(v_1)\otimes v_2$ on account of $\gamma$ being a $U$- comodule map. Further $\gamma(v_1)\cdot x=\epsilon(\gamma(v_1))x$ for every $x\in L(r)$ because both $E^{(\ell)}$ and $F^{(\ell)}$ annihilate $L(r)$ for every restricted $r$  (a rank 1 phenomenon). But then we see that $\rho(\gamma(v))(x\otimes y)=\epsilon(v_1)x\otimes v_2y=x\otimes vy$ as needed.

We turn now to the last assertion in (iii). Fix two primitive vectors $v^+$ and $w^+$ in $L(r)$ and $\overline{L}(\alpha)$, respectively. By the preceeding argument $v^+\otimes w^+$ is a primitive vector in $L(r)\otimes \overline{L}(\alpha)^{\text Fr}$. It also implies $K_{\bullet}(v^+\otimes w^+)=Kv^+\otimes w^+=\zeta^rv^+\otimes w^+$, where the ``$\bullet$" denotes the action of $U_{\zeta}$ in $L(r)\otimes \overline{L}(\alpha)^{\text Fr}$. As for the action of $\left[K\atop \ell \right]$, recall that $\left[K\atop \ell \right]= H - f$ in notation of lemma 1. Hence $\rho(H)=\left[\rho(E^{(\ell)},\rho(F^{(\ell)}\right]=1\otimes [e,f]=1\otimes h$, while $\rho(f)=f\otimes 1$. Thus $f_{\bullet}(v^+\otimes w^+)=fv^+\otimes w^+=0$, and therefore $H_{\bullet}(v^+\otimes w^+)=v^+\otimes hw^+=\alpha v^+\otimes w^+$. This shows that the weight of $v^+\otimes w^+$ is $\lambda=(r,\alpha)$ It remains to notice that by [Ja, 5.8.1] our module is irreducible.
\qed

\head \S 3 Decomposition of the $U$- module $\rz$\endhead

We need to review the structure of the PIMs (principal indecomposable modules) for 
$\r_{\zeta}$. For a start, we recall that by a theorem of Curtis- Lusztig [Lu 1, 2] (or directly) $\r_{\zeta}$ has $\ell$- simple modules of the form $L(r)|_{\r_{\zeta}}$ with $r$ $\ell$- restricted. We let $P(r)$ denote the projective cover of $L(r)$ viewed as a $\r_{\zeta}$- module. We prefer to construct these modules by exploiting the comodule theory of the finitary dual $(U_{\zeta})^0$ of $U_{\zeta}$. We refer to section 5 for a fuller discussion of this issue. Let $\Bbbk_{\zeta}[SL(2)]$ be the $\zeta$- deformation of the coordinate algebra of $SL_2(\Bbbk)$ ([CK 1]). It turns out that $(U_{\zeta})^0\simeq \Bbbk_{\zeta}[SL(2)]$ (cf. \S 5). Green's theory [Gr] guarantees existence of the injective hull $I_S$ for every right $\Bbbk_{\zeta}[SL(2)]$- comodule $S$. By [CK 1], or using the isomorphism just mentioned, we know that any simple comodule is isomorphic to $L(m)$ for some $m\in \Bbb{Z}^+$, treated as a right $\Bbbk_{\zeta}[SL(2)]$- comodule. Going in the opposite direction we can view $I(m)\colon =I_{L(m)}$ as a left $U_{\zeta}$- module. Restricting $m$ to the interval $0\le r< \ell$ we arrive at the $\ell$ $U_{\zeta}$- modules $I(r)$. Restricting  $I(r)$ to $\r_{\zeta}$ we obtain the $P(r)$. This can be seen as follows. Thanks to [APW 2, 4.6; Li, 6.3] the restriction $I(r)|_{\rz}$ is the injective hull of $L(r)|_{\rz}$. As $\rz$ is Frobenius, $I(r)|_{\rz}=P(r)$. An alternate construction of PIMs can be found in [Su].

To describe the structure of the $P(r)$ we need to bring in a new family of $U_{\zeta}$- modules. For an $m\in\Bbb{Z}^+$ we denote by $W(m)$ the $U_{\zeta}$- module generated by a primitive vector $v$ of weight $m$ subject to conditions $F^{(j)}v=0$ for $j> m$ and $F^{(i)}v\not= 0$ for all $i\le m$. We refer to $W(m)$ as the $m$th {\it Weyl module} for $U_{\zeta}$. The structure of $I(r)$ is as follows [CK 1]or [Ch 2]. Let $\rho$ be the reflection of $r$ in $\ell-1$, i.e. $\rho(r)=2(\ell-1)-r$. For every $r\not=\ell-1$, $I(r)$ is characterized as the unique extension of $W(\rho(r))$ by $L(r)$, and $I(\ell-1)=L(\ell-1)$. Further $I(r),\,r\not=\ell-1$, is a uniserial module with the factors $L(r),L(\rho(r)), L(r)$ in that order. 

Put $r'=\ell-2-r$. Then $\rho(r)=r'+\ell$, and by Steinberg- Lusztig theorem [Lu 1, 7.4] we have
$$ L(\rho(r))\simeq L(r')\otimes L(1)^{\text{Fr}}$$
Restricting to $\rz$ we get $L(\rho(r))|_{\rz}\simeq L(r')\oplus L(r')$. This yields in turn
$$P(r)\approx 2L(r)\oplus 2L(r')\tag 3-1$$
where $\approx$ signifies that two modules have the same composition factors. It follows readily that PIMs $P(r)$ and $P(s)$ are linked if and only if $s=r'$ or $s=r=\ell-1$. Therefore the block $\r_r$ of $\rz$ containing $P(r)$ has only one more PIM, namely, $P(r')$. Further, it is elementary that a PIM $P$ associated with the absolutely simple module $L$ has multiplicity $\dim L$ in its block. Thus $\r_r$ is the direct sum of $r+1$ copies of $P(r)$ and $r'+1$ copies of $P(r')$.

The next theorem requires an explicit decomposition of $\r_r$ into the direct sum of PIMs, and a special basis for $P(r)$.  

We start by describing a special basis for $P(r)$. Let $V$ be a nonsplit extension of $W(\rho(r))$ by $L(r)$. We recall that every finite dimensional $\Uz$- module is semi-simple as a $\U$- module [APW, \S 9]. 

Therefore $V|_{\U}\simeq L(r)|_{\U}\oplus W(\rho(r))|_{\U}$. 
It follows that $V$ contains an element $z_0$ of weight $r$ which is not in $W(\rho(r))$. The image $\overline{z}_0$ of $z_0$ in $L(r)$ is a primitive vector there, hence the set $\{F^{(i)}\overline{z}_0|i=0,\ldots,r\}$ forms a basis for $L(r)$. Note that $E\overline{z}_0=F\overline{z}_r=0$. Let's pull the $F^{(i)}\overline{z}_0$ into $V$ by setting $z_i=F^{(i)}z_0, i=0\ldots,r$. Then $Ez_0$ and $Fz_r$ lie in $W(\rho(r))$. Clearly their weights are $r+2$ and $-r-2$, respectively. Let $v_0$ be a primitive generator of $W(\rho(r))$, so that $\{v_i=F^{(i)}v_0|i=0,\ldots, \rho(r)\}$ is a basis of $W(\rho(r))$. Call such a basis {\it standard}. Then weight considerations lead to $Ez_0=av_{\ell-r-2}$ and $Fz_r=bv_{\ell}$ for some $a,b\in \Bbbk$. These are ``structure constants" of the extension. They depend on the choice of a standard basis for $W(\rho(r))$. We remark that neither $a$ nor $b$ is zero. For, the socle of $W(\rho(r))$ is the subspace spanned by $\{v_{r'+1},\ldots, v_{\ell-1}\}$, and this is also the socle of $V$ on account of $V$ being nonsplit extension. Now, were $a=0$, $z_0$ would be a primitive vector of $V$. Then $z_0$ would generate submodule of $V$ missing $\text{soc}\;V$, a contradiction. Likewise, $b=0$ implies the submodule $\Uz z_r$ misses $\text{soc}\;V$. We mention in passing that in the standard basis generated by $E^{(r'+1)}z_0$ the structure constants equal $[\ell-r-1]$ and $(-1)^r[\ell-r-1]$, respectively. It follows that $I(r)$ is a unique nonsplit extension of $W(\rho(r))$ by $L(r)$.

The desired basis is this. Define $\rho(r)$ vectors $w_k$ by the formulas
$$\align w_j=E^{(r'+1-j)}z_0, j=0,1,\ldots,r',\;w_{r'+1+i}&=v_{r'+1+i},\;i=0,1,\ldots,r\\  w_{\ell+i}=F^{(r+1+i)}z_0,\;i=0,1,\ldots,r'\tag 3-2\endalign$$ 
Thanks to the fact that $ab\not= 0$, the $\;w_j$ are nonzero scalar multiples of $v_j$, hence they form a basis of $P(r)$

We proceed to decomposition of $\r_r$. Let
$$B=\{0,1,\ldots, \frac{\ell-3}{2}\}\;\text{and}\quad\hat B=B\cup \{\ell-1\}$$
Recall the quantum Casimir element $c_{\zeta}$. It can be checked that $c_{\zeta}$ acts on $L(r)$ by multiplication by the scalar $\lambda_r=\frac{\zeta^{r+1}+\zeta^{-(r+1)}}{(\zeta-\zeta^{-1})^2}$. One can see easily that $\lambda_r=\lambda_s$ iff $r+s=\ell-2$. This limits the set of subscripts on $\lambda_r$ to $\hat B$. Let 
$$\Phi(x)=(x-\lambda_{\ell-1})\prod_{r\in B}(x-\lambda_r)^2$$ 
It turns out that $\Phi(x)$ is the minimal polynomial of $c_{\zeta}$ over $\Bbbk$. For, suppose that $f(x)$ is a polynomial annihilating $c_{\zeta}$. Then from $f(c_{\zeta})L(r)=0$ it follows that $f(\lambda_r)=0$. This being the case for all $r\in\hat B$, $f(x)$ is divisible by $\prod _{r\in\hat B}(x-\lambda_r)$. In the same vein, $f(c_{\zeta})P(r)=0$ implies that $f(x)$ is divisible by $\Phi(x)$. On the other hand $\Phi(c_{\zeta})=0$, because $\Phi(c_{\zeta})P(r)=0$ for all $r\in\hat B$.

It follows immediately that the algebra $\Bbbk[c_{\zeta}]$ is isomorphic to the direct sum of algebras $\Bbbk[x]/((x-\lambda_r)^2),\,r\in B$ and $\Bbbk$. Let $\epsilon_r,\,r\in B,\,\epsilon_{\ell-1}$ be idempotents corresponding to summands isomorphic to $ \Bbbk[x]/((x-\lambda_r)^2),\,r\in B$ and $\Bbbk$, respectively. The number of summands $\rz\epsilon_r$ of $\rz$ equals to the number of blocks, so that each $\rz\epsilon_r$ is block algebra. A simple verification yields that $\epsilon_r$ annihilates every indecomposable $P(s)$ with $s\not= r,r'$. Thus $\rz\epsilon_r=\r_r$ for every $r\in \hat B$.  

Recall the idempotents $e_j,\,j=0,1,\ldots,\ell-1$, defined in Remark 3. They give rise to the decomposition
$$\r_r=\oplus_{0\le j\le\ell -1}\r_re_j$$
In fact every $\r_re_j$ is a PIM for $\r_r$, because, as we mentioned earlier, $\r_r$ is a direct sum of $\ell$ PIMs. 

\proclaim{Theorem 3} (1) In the foregoing notation the $\rz\epsilon_re_j$ are $U$- stable.

(2) Let $\r_re_j$ be a summand of $\r_r$ isomorphic to $P(r)$. As a $U$- module $\r_re_j$ is a direct sum of $2(r+1)$ copies of the trivial representation and 2(r'+1) copies of the defining representation of $U$, i.e.,
$$\r_re_j\simeq \bar{L}(0)^{2(r+1)}\oplus \bar{L}(1)^{2(r'+1)}$$ 
\endproclaim
\demo{Proof} (1) Since $E^{(\ell)}$ and $F^{(\ell)}$ commute with $\epsilon_r$ and $e_j$, this statement is clear.

(2) We start with a PIM $P$. Let $N$ be the radical of $\rz$. Then $NP$ is the maximal submodule $W$ of $P$. Since $W$ is not semisimple, $N^2P$ is nonzero; hence it equals the socle of $P$

For the remainder of the proof set $P=\r_re_j$. By Lemma 3 \S2 $N^2$ is $U$- stable, therefore $N^2e_j$ is $U$- stable. We conclude that $\text{soc}\,P$ is $U$- stable. Next we note that as $E^{(\ell)}$ and $F^{(\ell)}$ commute with $K$ for every $\Uz$- module $M$, the $K$- weight subspace $M_\lambda,\,\lambda\in \Bbb Z_{\ell}$, is stable under $U$. Take $M=\text{soc}\,P$. Then $M_{\lambda}$ is $1$- dimensional or zero. Therefore every vector of $\text{soc}\,P$ is $U$- trivial, i.e. they are fixed points under the action of $U$. Moreover, the formula (3-1) shows that, for every weight $\lambda$, $P_{\lambda}$ is $2$- dimensional. Since the bottom composition factor has only one equal composition factor we see that the elements of $P_{\lambda}$, $\lambda$ a weight of $\text{soc}\,P$ are fixed by $U$. As $\text{soc}\,P\simeq L(r)$ we get $r+1$ values of $\lambda$. This accounts for $2(r+1)$ copies of $\bar{L}(0)$ in $\r_re_j$.

Suppose $\r_re_j\simeq P(r)$ and recall the basis $w_j,\,j=0,1\ldots, r'$ as in (3-2) above. It remains to show that none of the $w_j$ are $U$- trivial. Let $z_0$ be a generator of weight $r$ of $\r_r$. Then $w_j=E^{(p)}z_0$ for the $p$ such that $j+p=r'+1$. Keeping in mind that $z_0$ is a fixed point we have
$$[F^{(\ell)},w_j]=[F^{(\ell)},E^{(p)}]z_0=-\left(\sum_{0\le i\le p} F^{(\ell-i)}\left[{K;2i-\ell-p}\atop i\right]E^{(p-i)}\right)z_0$$
We claim that all terms in the right- hand side of that formula with $i<p$ are zero. 

For, suppose $p-i>0$. Up to a nonzero scalar $F^{(\ell-i)}E^{(p-i)}$ equals $F^{(\ell-p)}F^{(p-i)}E^{(p-i)}$. Further, by formula (1-1) $F^{(p-i)}E^{(p-i)}z_0=g(c_{\zeta}, K)(c_{\zeta}-\lambda_r)z_0$ for a polynomial $g(x,y)$. Also it is immediate that $(c_{\zeta}-\lambda_r)z_0$ lies in $\text{soc}\,\r_re_j$. Now $F^{(\ell-p)}$ annihilates $\text{soc}\,\r_re_j$ because $\ell-p>r$ for every $p$.

For the $i=p$ term we have  
$$\align -F^{(\ell-p)}\left[{K;p-\ell}\atop p\right]z_0=-\left[{r+p-\ell}\atop p\right]F^{(r+1+j)}z_0&=-\left[{-(j+1)}\atop p\right]w_{\ell+j}\\&=(-1)^{p+1}\left[{r'+1}\atop p\right]w_{\ell+j}\endalign$$
As $\left[{r'+1}\atop p\right]\ne 0$ we are done
\qed
\enddemo

As a by-product of the proof we note

\proclaim{Corollary} (1)$N^2$ is generated by all $(c_{\zeta}-\lambda_r)\epsilon_r,\,r\in B$

(2) $N^2$ consists of $U$- trivial elements
\endproclaim
\demo{Proof} Both statements follow from the fact that $N^2$ is the sum of $\text{soc}P$, over PIMs $P$ , and the fact that $\text{soc}P(r)$ is generated by by $(c_{\zeta}-\lambda_r)z_0$, where $z_0$ is a weight $r$ generator of $P(r)$
\qed
\enddemo

\head \S 4 The Center of $\Uz$\endhead

We let by $Z(A)$ denote the center of a $\Bbbk$- algebra $A$. We can use central idempotents $\epsilon_r$ to split $\Uz$ into the direct sum of subalgebras $\epsilon_r\Uz$. This leads to the splitting of $Z(\Uz)$ into the direct sum of $Z(\epsilon_r\Uz)$, which reduces the description of $Z(\Uz)$ to $Z(\epsilon_r\Uz)$.

The smash product representation of $\Uz$ implies that $\epsilon_r\Uz$ is isomorphic to $\r_r\# U$ where $\r_r=\epsilon_r\rz$ is a block of $\rz$. When $r=\ell-1$, $\r_{\ell-1}=\text{End}\, L(\ell-1)\simeq M_{\ell}(\Bbbk)$. Let $N$ be the radical of $\rz$ as before. Denote by $c=4fe+(h+1)^2$ the usual Casimir element of $U$. Since $\r_{\ell-1}\cap N=0$ we derive from Lemma 3 that $\epsilon_{\ell-1}\Uz\simeq M_{\ell}(\Bbbk)\otimes U\simeq M_{\ell}(U)$. It follows that $Z(\epsilon_{\ell-1}\Uz)=1\otimes Z(U)=1\otimes \Bbbk[c]$.

In what follows $r\ne \ell-1$. We need to review a description of the center of $\r_r$. The result may be well-known, but we don't have a specific reference. The statement is that $Z(\r_r)$ is three- dimensional. We sketch a proof of it and also construct a basis of the center.

Let $\frak b$ be the basic algebra of $\r_r$. As in ([Su]) $\frak b$ is $\Bbbk$- algebra on a basis 
$$\{e,c,a,b,e',c',a',b'\}$$
subject to relations
$$\align ec&=ce=c\\
         ab'&=ba'=c\\
         e'c'&=c'e'=c'\\
         b'a&=a'b=c',\endalign$$ 
$e$ is a left unit of all non-primed generators and a right unit of all primed generators, $e'$ is a left unit of all primed generatora and a right unit of all non-primed generators, and all other products of basic elements are zero.

It is immediate from the definition of $\frak b$ that $\{1_{\frak b},c,c'\}$ spans $Z(\frak b)$. On the other hand it is straightforward to check that 
$$\epsilon_r,\epsilon_r(c_{\zeta}-\lambda_r),\,\text{and}\;\theta=\epsilon_r(c_{\zeta}-\lambda_r)(\sum_{j=0}^{r}e_{r-2j})\tag 4-1 $$
lie in the center of $Z(\r_r)$ and are linearly independent. We conclude that the set (4-1) is a basis of $Z(\r_r)$. It follows readily that $M=\Bbbk (c_\zeta-\lambda_r)\epsilon_r\oplus\Bbbk\theta$ is the unique maximal ideal of $Z(\r_r)$

We proceed to the proof of

\proclaim{Theorem 4} $Z(\epsilon_r U)={\Bbbk}\epsilon_r\oplus (M\otimes {\Bbbk }[c])$
\endproclaim
\demo{Proof} In the easy direction we want to show that $(c_\zeta-\lambda_r)\epsilon_r\otimes {\text k}[c]$ and $\theta\otimes {\text k}[c]$ lie in $Z(\epsilon U)$. To this end we note that for every $x\in N^2,\,y\in U$, and $w\in\r_r$ we have $(xy)w=xwy+x[y,w]=xwy$, because by Lemma 3 $[y,w]\in N$ and $N^3=0$. Further every $x\in Z(\r_r)$ commutes with every element of $U$, because it lies in $\Bbbk[c_{\zeta},K]$. By Corollary of Theorem 3 the assertion follows.

We now turn to the harder part of the proof. Let $z=\sum x_iy_i,\;x_i\in {\bold u}_r,\;y_i\in U$ be a central element. We may assume that both sets $\{x_i\}$ and $\{y_i\}$ are linearly independent over $\Bbbk$. Then, as the $y_i$ commute with $K$, $Kz=zK$ implies $Kx_i=x_iK$ for every $i$. But every element of $\rz$ that commutes with $K$ lies in the subalgebra of $\rz$ generated by $c_\zeta$ and $K$. Further, this subalgebra centralizes $U$. It follows that the $y_i$ lie in the center of $U$.

For the remainder of the proof we abbreviate $E^{(\ell)}$ to $e$, $F^{(\ell)}$ to $f$ and we write $h$ for $H$. 
Let $V=\Bbbk t\oplus \Bbbk v$ be a typical two-dimensional $U$- submodule of $\r_{\zeta}$ as in Theorem 3. Then $V\# U$ is a $U-U$ subbimodule of $U_{\zeta}$. Consequently, the left multiplication $\lambda_u$ by an $u\in U$ induces a right $U$- linear map in $V\# U$. We denote by $M_u\in M_2(U)$ the matrix of $\lambda_u$ relative to the basis $\{t,v\}$. We write $M$ for $M_c$, hence $M^n$ for $M_{c^n}$. We put $M^n=\pmatrix a_n&b_n\\c_n&d_n\endpmatrix$ and denote by $\omega$ the standard automorphism of $\frak{sl}_2$ defined by $\omega(e)=f,\,\omega(f)=e,\,\omega(h)=-h$. To proceed with the proof, we need
 
\proclaim{Lemma 4} (1) The elements of $M^n$ are determined by $a_n$. We have $b_n=[f,a_n],\,c_n=\omega(b_n),\,d_n=\omega(a_n)$.

(2) $a_n$ is a polynomial in $c$ and $h$. Namely, $a_n=\displaystyle\sum_{i=0}^nc^i\phi_i(h)$ with $\phi_i(h)$ of degree $\le 1$. Further, $\phi_n(h)=1$ and $\phi_{n-1}(h)=2nh+n(2n+1)$.

\endproclaim
\demo{Proof}For calculations below one must keep in mind the action of $e,f,h$ on $t$ and $v$ as well as the definition of $c$. The action is given by $[e,t]=0,\;[f,t]=v,\;[h,t]=t$ and $[e,v]=t,\;[f,v]=0,\;[h,v]=-v$. By definition we have $(\ast)\;c^nt=ta_n+vc_n$. We apply $\text{ad}_L\ h$ and $\text{ad}_L\ f$ in turn to both sides of $(\ast)$. For $\text{ad}_L h$ we get  $c^nt=t(a_n+[h,a_n])+vg$ for some $g\in U$, which implies $[h,a_n]=0$ on account of $t,v$ being linearly independent over $U$. Thus $a_n$ lies in the subalgebra of $U$ generated by $c$ and $h$, as claimed in the first part of (2). Further $[f,c^nt]=c^n[f,t]=c^nv$. On the other hand, $[f,ta_n]+[f,vc_n]=t[f,a_n]+vg$ for a $g\in U$. It follows that $b_n=[f,a_n]$. As for the last claim in (1), it suffices to do the $n=1$ case.

This can be computed directly. The result is 
$$ct=t(c+2h+3)+4ve\;\text{and}\;cv=4tf+v(c-2h+3)$$ confirming the last part of (1). 

The remaining assertions in (2) hold for $n=1$ by the preceeding two equations. The general case can be verified by induction on $n$. We leave the details to the reader
\qed
\enddemo

Conclusion of the proof of the Theorem. Let $z$ be an element of $Z(\epsilon U)$. From the previous remarks we have 
$$z=\sum\limits_{i\le n}x_i\cdot c^i,\;x_i\in \Bbbk[c_{\zeta}, K]$$
The theorem's statement is equivalent to the claim that for every $i\ge 1$

$(*)$ $x_i$ lies in the span of $(c_\zeta-\lambda_r)\epsilon$ and $\theta$

By the first part of the proof we know that if $x_i$ satisfies $(*)$ then $x_ic^i$ is in the center of $U$. Hence we may assume that the leading coefficient $x_n$ of $z$ doesn't satisfy $(*)$. Next we compute $tz=zt$ in two ways. First, $tz=\sum tx_ic^i$. Second, by Lemma 4 $$zt=x_ntc^n+[x_nt(2nh+n(2n+1)+x_nv(4ne)+x_{n-1}t]c^{n-1}+\text{lower $c$- degree terms }$$. 
Since the powers $c^i$ are independent over $\r_r\otimes\Bbbk[h]$ we can equate the coefficients of $c^n$ and $c^{n-1}$ on both sides of the equation obtaining
$$\align\text {(i)}\quad tx_n&=x_nt\\\text{(ii)}\;tx_{n-1}&=x_nt(2nh+(2n+1)n)+x_nv(4ne)+x_{n-1}t\endalign $$
Now (i) implies that $x_n\in Z({\bold u}_r)$. For, it shows that $x_n$ commutes with $t$. Replacing $t$ with $v$ and repeating the argument above we conclude that $x_n$ commutes with $v$. On the other hand if $w$ is $U$- trivial, then $zw=wz$ clearly implies $x_nw=wx_n$. Since by Theorem 3 $\r_r$ has a basis whose elements fall into either of these two types we get the assertion.

(ii) shows that $x_nt=x_nv=0$ as we are working in the tensor product $\r_r\otimes U$. Also by (i) we can write $x_n=a_0\epsilon+a_1 (c_\zeta-\lambda_r)\epsilon+a_2\theta$. Taking into account $(c_\zeta-\lambda_r)t=\theta t = 0$, because both $(c_\zeta-\lambda_r)$ and $\theta$ lie in $N^2$ while $t,v\in N$ and $N^3=0$, we deduce $x_nt=a_0\epsilon t=a_0t=0$. This forces $a_0=0$, a contradiction 
\qed
\enddemo

We proceed to consider the block decomposition of $\Uz$. Let $\epsilon$ be a primitive central idempotent of $\rz$. Then $\epsilon \Uz$ is a two-sided ideal of $\Uz$. In fact, it is a block of $\Uz$ in the sense of the following

\proclaim{Corollary} $\epsilon \Uz$ is an indecomposable algebra
\endproclaim
\demo{Proof} By the theorem the center of $\epsilon\Uz$ is a local algebra
\qed
\enddemo

\head \S 5 Cofinite Ideals of $U$\endhead

In this section we change our previous notation in that we are going to write $U$ for $\Uz$ and $H$ for $\Bbbk_{\zeta}[SL(2)]$. $U^0$ now stands for the finitary dual of $U$. The starting point of the section is existence of a Hopf pairing ([Ta 1]) between $U$ and $H$. An equivalent formulation is existence of a Hopf map $\psi\colon H\to U^0$. The construction of $\psi$ runs as follows (cf. [Ta 2], [deCL]). Let $L(1)$ be the $2$- dimensional representation of $U$. It gives rise to the coordinate functions $c_{ij},\,i,j=1,2$. Let $A$ be the subalgebra of $U^0$ generated by the $c_{ij}$. Then $A$ is a Hopf subalgebra of $U^0$ and the natural mapping $x_{ij}\to c_{ij},\,i,j=1,2$, where $x_{ij}$ are the standard generators of $H$ [CK 1], extends to the Hopf algebra map $H\to A$. In fact, that map is an isomorphism \linebreak [Ta 1], and moreover according to [APW 1, Appendix] we have $A=U^0$. Below we will identify $H$ with $A$

Let $\pi$ be the natural algebra homomorphism $U\to (U^0)^{\ast}$ [Sw, 6.0], $\pi(u)(\alpha)=\alpha(u),\,u\in U,\alpha\in U^0$. Note that the image of $\pi$ is in $(U^0)^0$.
 
 We compose $\pi$ with $\psi^0\colon (U^0)^0\to H^0$ to obtain $\phi=\psi^0\pi\colon U\to H^0$. Explicitly, $\phi(u)(x)=((\psi^0\pi)(u))(x)=\pi(u)(\psi(x))=(\psi(x))(u)$. 

Using $\psi$ we define a bilinear form $$\langle\,,\,\rangle\colon U\otimes H\to\Bbbk,\,\langle u,x\rangle=\psi(x)(u)$$ 
For a subspace $L$ of $U$ we set $L^{\perp}=\{x\in H|\langle L,x\rangle=0\}$, and likewise for an $M\subset H$. We refer to the above as the annihilators of $L$ and $M$, respectively. We recall that a subspace $F$ of $H^{\ast}$ is said to be dense in $H^{\ast}$ if $F^{\perp}=0$ in the natural pairing $H\times H^{\ast}\to\Bbbk$. A subspace $L$ of $U$ ($M$ of $H$) is said to be closed if $L^{\perp\perp}=L\,(M^{\perp\perp}=M)$. A pairing $U\times H^{\ast}\to\Bbbk$ is nondegenerate on the left (right) if $H^{\perp}=0$ $(U^{\perp}=0)$. A pairing is nondegenerate if it is both left and right nondegenerate.

\proclaim{Lemma 5} Every cofinite ideal of $U$ and every finite- dimensional subspace of $H$ is closed.
\endproclaim
\demo{Proof} Let $I^{\perp\ast}$ stand for the annihilator of $I$ in $U^{\ast}$. In view of $H=U^0$ we have $I^{\perp\ast}=I^{\perp}$. But then $I^{\perp\perp}=I^{\perp\ast\perp}$ which is $I$ by [Sw, A.1].

The second statement is well-known for a nondegenerate pairing [Ab, 2.2]. The same proof works for a right nondegenerate pairing. It remains to note that our pairing is right nondegenerate by the equality $H=C$. 
\qed
\enddemo

We let $\Lambda_{\text cof}(U)$ and $\Lambda^{\text fin}(H)$ denote the lattices of cofinite ideals of $U$ and finite-dimensional subcoalgebras of $H$, respectively.

\proclaim{Lemma 6} (1) For every ideal $I$ of $U$, and every subcoalgebra $C$ of $H$ $I^{\perp}$ is a subcoalgebra of $H$ and $C^{\perp}$ is an ideal of $U$

(2) $\Lambda_{\text cof}(U)$ is antiisomorphic with $\Lambda^{\text fin}(H)$ under $I\mapsto I^{\perp}$
\endproclaim
\demo{Proof} (1) The usual proof ([Sw, 1.4.3]) of the assertions works here thanks to the property $\phi(U)$ dense in $H$, which is equivalent to $U^{\perp}=0$

(2) It is straightforward to see that the mappings $I\mapsto I^{\perp}$ and $C\mapsto C^{\perp}$ set up two inclusion reversing correspondences between $\Lambda_{\text cof}(U)$ and $\Lambda^{\text fin}(H)$. Lemma~5 makes it clear that those maps are mutual inverses of each other. By general principles they send unions to joins and conversely.
\qed
\enddemo

For a coalgebra $C$ we let $C^C$ denote the right regular $C$- comodule. We put $E=\text{End}^C(C^C)$ for the algebra of all $C$- endomorphisms of $C^C$. We recall that $C$ has the natural structure of the right $C^{\ast}$- module via $c\leftharpoonup c^{\ast}=\sum\langle c_1,c^{\ast}\rangle c_2$. 

The next very useful lemma is a variation on a well-known statement in the theory of coalgebras. A weaker form can be found in [Ho, p.18]. 

\proclaim{Lemma 7} $E$ is antiisomorphic with the algebra $C^{\ast}$ acting on $C$ by the right `hits'
\endproclaim
\demo{Proof} In one direction, given $f\in E$ we associate with it $\phi=\epsilon\circ f\in C^{\ast}$. In the opposite direction, we send $\phi\in C^{\ast}$ to $f_{\phi}\colon x\mapsto x\leftharpoonup \phi$. One can check that $f_{\phi}\in E$ and the above defined maps are mutually inverse antiisomorphisms between $E$ and the image of $C^{\ast}$ in $E$. It remains to note that by [Sw, 9.1.2] $\phi\mapsto f_{\phi}$ is an injection.
\qed
\enddemo

A description of the lattice $\Lambda^{\text fin}(H)$ begins with  the representation of $H$ as a direct sum of injective indecomposables. We view $H$ as a coalgebra only. We know that every simple right $H$- comodule is $L(r)$ for some $r$ by [CK 1] (or using the fact that $H=U^0$ and the classification of the simple finite-dimensional $\Uz$- modules. Let $I(r)$ be the injective hull of $L(r)$, and put $m(r)=\dim L(r)$. Then we have by [Gr] ( keeping in mind that the $L(r)$ are absolutely simple)
$$H=\bigoplus I(r)^{m(r)}\tag 5-1$$

We refer to the representation (5-1) for a coalgebra $C$ as {\it the Green decomposition} for $C$.

We denote by $I(r)^{(j)},\,j=1,\ldots,m(r)$, the $j$-th copy of $I(r)$ in $H$. Let $C$ be a subcoalgebra of $H$. We want to connect the structure of a subcoalgebra to that of injective comodules. The first step is

\proclaim{Proposition 2} (i) $C=\bigoplus_r(\bigoplus_j C\cap I(r)^{(j)})$

(ii) The set of subcomodules $C\cap I(r)^{(1)}$ determines $C$
\endproclaim
\demo{Proof} (i) holds iff every component of an $x\in C$ lies in $C$. Hence assume $x=\sum x_{r,j}$ and let $\pi_{r,j}$ be the projection of $H$ on $I(r)^{(j)}$. By lemma 7 we can regard $\pi_{r,j}$ as an element of $H^*$. As $C\leftharpoonup \pi_{r,j}\subset C$, the proof is complete.

(ii) It suffices to show that for every isomorphism $\phi\colon I(r)^{(1)}\simeq I(r)^{(j)}$,\   $(C\cap I(r)^{(1)})\phi=C\cap I(r)^{(j)}$. Since $\phi$ can be lifted to an $H$- endomorphism of $H$ we can view $\phi$ as an element of $H^*$. But then $(C\cap I(r)^{(1)})\phi\subset C\cap I(r)^{(j)}$ for the same reason as above. Since $\phi$ has an inverse $I(r)^{(j)}\to I(r)^{(1)}$, the equality follows
\qed
\enddemo

We can improve significantly on the above proposition. We need to introduce some notation and terminology. Given a right $H$- comodule $M$ we say that a subcomodule $N$ of $M$ is {\it fully invariant} if it is stable under the action of $\text{End}^H(M)$. In the case $M=H$ let denote by ${}^HN$ and $C(N)$ the left $H$-comodule and the subcoalgebra generated by $N$, respectively. Of course, $C(N)={}^HN$. We let $\Lambda_r$ denote the sublattice of the fully invariant subcomodules of $I(r)$. Let $\prod\Lambda_r$ be the Cartesian product of the $\Lambda_r$ and put $\varphi_r$ for the natural projection $\prod\Lambda_r\to \Lambda_r$. Set $E_{r,s}=\text{Hom}^H(I(r),I(s))$ and call an element $(X_r)_{r\in\Bbb Z^+}\in \prod\Lambda_r$ {\it balanced} if 
$$X_rE_{r,s}\subset X_s\,\text{for all}\;r,s$$
Let $\Lambda^{\text{coalg}}(H)$ be the lattice of the subcoalgebras of $H$. Also, in what follows we drop the superscript and write $I(r)$ for $I(r)^{(1)}$

\proclaim{Proposition 3} The mapping $\theta\colon C\mapsto (C\cap I(r))_{r\in\Bbb Z^+}$ is a lattice embedding $\Lambda^{\text{coalg}}(H)\to \prod\Lambda_r$. The image of $\theta$ is the set of all balanced elements of $\prod\Lambda_r$ and $\varphi_r\theta$ is a surjection for all $r$
\endproclaim
\demo{Proof} As in the proof of Proposition 2, for every $\phi\in E_{r,s}$ we have $(C\cap I(r))\phi\subset C\cap I(s)$. Thus $C\cap I(r)\in \Lambda_r$ for every $r$ and also the sequence $\theta(C)$ is balanced.

The rest of the proposition is based on several observations. Let $\pi_r$ be the projection of $S$ on $I(r)$. Then the following holds (i) $C\cap I(r)=C\pi_r$,\ (ii) For every $X\in \Lambda_r,\;C(X)\pi_r=X$, and (iii) $C(X)\pi_s=XE_{r,s}$

(i) is obvious. As for (ii) it suffices to show that $C(X)\pi_r\subset X$. We have $C(X)\pi_r=({}^SX)\pi_r=XS^*\pi_r$. Pick a $\phi\in S^*$. Then $\phi\pi_r$ maps $X$ into $I(r)$. By injectivity of $I(r)$ $\phi\pi_r$ is the restriction of a $\psi\in \text{End}(I(r))$. As $X$ is fully invariant, (ii) follows.

For (iii) we have $C(X)\pi_s=XS^*\pi_s$. As $S^*\pi_s|_{I(r)}=E_{r,s}$, we are done.

Now, pick two subcoalgebras $C$ and $D$. For every $r\in\Bbb Z^+$ we have $(C+D)\cap I(r)=(C+D)\pi_r=C\pi_r+D\pi_r=C\cap I(r)+D\cap I(r)$ which proves that $\theta$ is a lattice map.

Further, suppose $(X_r)_{r\in\Bbb Z^+}$ is a balanced sequence. Put $C=\sum C(X_r)$. Then $C\pi_r=\sum C(X_s)\pi_r=X$ on account of (ii) and (iii). Also, (ii) shows that $\varphi_r\theta$ is a surjection
\qed
\enddemo

The last proposition points to the importance of knowing fully invariant subcomodules of $I(r)$. It will turn out that these coincide with the subcomodules of $I(r)$. Thus we proceed to a description of the subcomodule lattice of $I(r)$. In this we follow [CK 1, Thm. 5.2], also see [Ch 2]

There is a symmetry at work here expressed, in general, by the $\ell$- affine Weyl group (cf. [Th]). In the case at hand we will use instead a bijection $\rho\colon \Bbb Z\to\Bbb Z$, called an $\ell$- reflection. Write $m=m_0+m_1\ell$, where $0\le m_0\le \ell-1$. Define $\rho(m)=m$, if $m_0=\ell-1$ (i.e. $m$ is a Steinberg weight). If $m_0\not= \ell-1$, then set $\rho(m)$ equal to the reflection of $m$ in the nearest Steinberg weight to the left. All in all we have $$\rho(m)=\ell-2-m_0+(m_1-1)\ell\;\text{if}\;m_0\not=\ell-1,\,\text{and}\,\rho(m)=m\;\text{otherwise}$$

Next we need to bring in a new family of $S$- comodules $\{M(r)|r\in\Bbb Z^+\}$. By definition $M(r)$ is the (contragredient) dual of $W(r)$, i.e., $M(r)=W(r)^*$. We note an independent construction of this family. Let $A=\Bbbk_{\zeta}[e_1,e_2]$ be the $\zeta$- plane, that is the algebra generated by $e_1,\,e_2$ subject to one relation $e_2e_1=\zeta e_1e_2$. $A$ is a ${\Bbb Z}^+$-graded algebra with the $n$-th component $A_n$ consisting of all homogeneous polynomials of degree $n$. $A$ is also a natural right $S$- comodule algebra via $$\omega(e_j)=\sum e_k\otimes x_{kj}$$
Each $A_n$ is a subcomodule and the family $\{A_r\}$ coincides with $\{M(r)\}$.

We now state the main facts about $I(r)$ from [CK 1] (see also [Ch 2]). 

\proclaim{Proposition 4} (1) $I(r)$ is a local, self-dual comodule. 

(2) If $r=\ell-1+r_1\ell$ is 
a Steinberg weight, then $I(r)=L(r)$. 

(3) If $r=r_0+r_1\ell$ with $r_0\ne\ell-1$ and $r_1>0$, then nonzero proper subcomodules of $I(r)$ are precisely $M(r)$, $W(\rho^{-1}(r))$, $L(r)$. 

(4) If $r=r_0$, the nonzero proper subcomodules of $I(r)$ are precisely $W(\rho^{-1}(r))$ or $L(r)$.
\endproclaim

We proceed to a description of ``links" between the $I(r)$.

\proclaim{Lemma 8} $\dim\text{Hom}^H(I(r),I(s))\le 1$ for all $r,s$. Equality holds iff $r=\rho(s)$ or $s=\rho(r)$.
\endproclaim
\demo{Proof} It is well-known  [Lu 1] that $W(m)$ has two composition factors $L(m)$ and $L(\rho(m))$. From Proposition 4 we see that $I(r)$ and $I(s)$ have a common factor iff $r$ and $s$ are as stated in the lemma.

Suppose $r=\rho(s)$ for definiteness. The preceeding proposition makes it clear that the multiplicity of $L(s)$ in $I(r)$ is $1$. Hence the dimension in question is at most $1$. Now dualize $0\to W(s)\to I(r)$ to get $I(r)^*\to W(s)^*\to 0$. As $I(r)$ is self-dual, we are done.
\qed
\enddemo

\proclaim{Lemma 9} Let $r=\rho(s)$ and $f_r, g_s$ be a non-zero elements of $\text{Hom}^H(I(r),I(s))$ and $\text{Hom}^H(I(s),I(r))$, respectively. Then $\text{Ker}\,f_r=M(r)$ and $\text{Ker}\,g_s=W(\rho^{-1}(s))$
\endproclaim
\demo{Proof} We only consider $f_r$. From the previous lemma the image of $f_r$ is $M(s)$. Since $\text{soc}\,M(s)=L(s)$ we have $f_r(M(r))=0$. As the codimension of $M(r)$ equals $\dim M(s)$ the assertion follows
\qed
\enddemo

We combine propositions 3 and 4 with the last lemma to derive the structure of an arbitrary finite-dimensional subcoalgebra of $H$.

We say that a coalgebra is local if it has a unique maximal subcoalgebra. Similarly, a comodule is local if it has a unique maximal subcomodule.

\proclaim{Theorem 5}The following properties hold for $\Lambda^{\text{fin}}(H)$
\roster
\item The lattice $\Lambda^{\text{fin}}(H)$ is distributive.

\item Every subcoalgebra is a unique sum of local subcoalgebras.

\item Every local subcoalgebra has the form $C(X)$, where $X$ is a local\linebreak subcomodule of an $I(r)$.

\item For a coalgebra $C=C(X)$, where $X$ is a subcomodule of $I(r)$, the Green decomposition is
$$C=X_{\rho(r)}^{m_{\rho(r)}}\oplus X_r^{m_r}\oplus X_{\rho^{-1}(r)}^{m_{\rho^{-1}(r)}}$$
with $X_{\rho^{\pm 1}(r)}=X_r\text{Hom}^H(I(r), I(\rho^{\pm 1}(r))$
\endroster
\endproclaim
\demo{Proof} (1) Proposition 4 makes it visible that the lattice of subcomodules of every $I(r)$ is distributive. Then, by a result in [Ste], we conclude that every subcomodule is fully invariant. As the class of distributive lattices is closed under taking of Cartesian products and sublattices, $\Lambda^{\text{fin}}(H)$ is distributive.

(2) This is a standard fact in the theory of distributive lattices \cite{DP}

(3) In one direction, assume $C=C(X)$ for a local $X\subset I(r)$. Were $C=D+E$ for some proper subcoalgebras $D$ and $E$ we would have by proposition 3, $$X=C\cap I(r)=(D+E)\cap I(r)=D\cap I(r)+E\cap I(r)$$,
a contradiction.

Conversely, pick a local $C$. Let $X_j=C\cap I(j)$. Then $C$ and $\sum C(X_j)$ have equal intersection with every $I(r)$ by Proposition 3, hence they are equal by Proposition 2. Thus $C=C(X_j)$ for some $X_j$. This forces $X_j$ to be local, for were $X_j=X_j'+X_j''$ with proper $X_j',\,X_j''$, so would $C$ be a sum of proper subcoalgebras.

(4) As $L(r)$ is the socle of $I(r)^{(j)}$ $C\cap I(r)^{(j)}$ is an indecomposable subcomodule of $C$. Proposition 2(i) gives $C$- injectivity of these and the Green decomposition of $C$. Further, properties (i)-(iii) in the proof of Proposition 3 say that\linebreak $C(X)=\displaystyle{\oplus\sb s}(XE_{r,s})^{m_s}$. But by lemma 8 $E_{r,s}\not= 0$ iff $s=\rho(r)$ or $s=\rho^{-1}(r)$
\qed
\enddemo

We pass on to a classification of the cofinite ideals. A few preliminaries are in order. 
Given an $H$- comodule $M$, $\text{ann}_UM$ is the annihilator of $M$ in $U$, where $M$ is regarded as a $U$- module via the duality $\langle,\rangle$. For an ideal $N$ of $U$ and a $U$- module $I$, $\text{ann}_IN$ denotes the submodule $\{x\in I|Nx=0\}$.

For an $H$- comodule $M$, $\text{cf}(M)$ denotes the {\it coefficient space} of $M$ defined in [Gr].

\proclaim{Lemma 10} If $M$ is a subcomodule of $H$, then $\text{cf}(M)=C(M)$, and  $C(M)^{\perp}=\text{ann}_UM$
\endproclaim
\demo{Proof} By the definition of $\text{cf}(M)$ one has $\text{cf}(M)^{\perp}=\text{ann}_UM$.

By ([Gr, 1.2f]) $\text{cf}(M)$ is a subcoalgebra containing $M$, whence $C(M)\subset \text{cf}(M)$. Conversely, pick $m\in M$ and let $\Delta(m)=\sum m'_i\otimes m''_i$. Then $m\leftharpoonup m_i^*=m_i''\in C(M)$ for every $i$. But, by definition $\text{cf}(M)$ is spanned by the various $m_i''$ as $m$ runs over a basis for $M$
\qed
\enddemo

We exploit the duality between $H$ and $U$ to obtain

\proclaim{Theorem 6} Let $X$ run over the set of injective indecomposable $H$- comodules and their local subcomodules.

Any cofinite ideal of $U$ is a unique intersection of the ideals $\text{ann}_UX$\qed
\endproclaim

Our next goal is to give the structure of the ideals $\text{ann}_UX$ of the previous theorem.

Two weights $r$ and $s$ are called {\it linked} if $r=\rho(s)$ or $r=\rho^{-1}(s)$. The linkage generates an equivalence relation on $\Bbb Z^+\times\Bbb Z^+$. The class of $r$, $\text{bl}(r)$, is called the block of $r$. If $r$ is not a Steinberg weight, i.e. $r_0\ne {\ell -1}$, then $\text{bl}(r)=\text{bl}(m)$, where $m$ is the smallest weight in $\text{bl}(r)$. Such $m$ is necessarily restricted, that is, $0\le m\le\ell-1$. The set $\text{bl}(r)$ equals $\{\rho^{-i}(r_0)|i\in\Bbb Z^+\}$. We remark that $\oplus_{s\in\text{bl}(r)}I(s)^{d_s}$ is exactly the block subcoalgebra summand of $H$ containing $I(r)$, whence our notation. We may ignore all injectives not in the block of $I(r)$. For those in the block we write $r_j=\rho^{-j}(r_0)$ and $I(j)=I(r_j)$. 

\proclaim{Lemma 11} Let $N$ be a cofinite ideal of $U$, $C=N^{\perp}$ and $I$ one of the $I(j)$. Put $X=C\cap I$. Then $X=\text{ann}_IN$
\endproclaim
\demo{Proof} For a subcomodule $Y\subset I$ $NY=0$ iff $N\subset \text{cf}(Y)^{\perp}$. Therefore $C=N^{\perp}\supset \text{cf}(Y)^{\perp\perp}=\text{cf}(Y)$. However $\text{cf}(Y)=C(Y)$ and therefore $C\supset C(Y)$ which implies $X\supset Y$ by Proposition 3. Thus $X$ is the largest subcomodule of $I$ annihilated by $N$
\qed
\enddemo

Let $X$ be as in Theorem 6 and let $$ X=X_1\supset X_2\supset\cdots X_t\supset 0\quad (t\le 4)\tag 5-2 $$ be a composition series for $X$. Let $\{L_i\}$ be the set of composition factors of (5-2) in order from top to bottom. Put $C_i=\text{cf}(L_i)$ and $\frak {m}_i=\text{ann}_UL_i$.

Recall that $A\wedge B$ denotes the wedge [Sw,9.0] of two subspaces $A$ and $B$.

\proclaim{Theorem 7} In the foregoing notation
 $$ C(X)=C_t\wedge \cdots\wedge C_1\;\text{or\ equivalently}\;\text{ann}_UX=\frak m_t\cdots m_1$$
\endproclaim
\demo{Proof} Recall that $\frak{m_i}=C_i^{\perp}$ by Lemma 10, while from [Sw, 9.0.0(b)], which holds for our pairing as well, we derive $(C_t\wedge \cdots\wedge C_1)^{\perp}=\frak m_t\cdots \frak m_1$.

Let $N=\frak m_t\cdots \frak m_1$ and put $D=N^{\perp}$. We will prove the equality $D=C(X)$ by showing that $D\cap I(k)=C(X)\cap I(k)$ for all $k$ (see Proposition 2). The original Hopf pairing $U\times H\to\Bbbk$ induces the pairing $U/N\times D\to \Bbbk$ which is clearly nondegenerate. It follows that the algebras $D^*$ and $U/N$ are isomorphic. As $\cap\frak{m}_i$ is the radical of $U/N$, we conclude that $$\text{corad}(D)=(\text{rad}\,D^*)^{\perp}=\sum\frak{m}_i^{\perp}=\sum C_i$$

Suppose $X\subseteq I(j)$. Apart from the trivial case $X=L(j)$, $X$ can be equal to $M(j),W(j+1)$ or $I(j)$. The first two cases are similar. Thus we consider two cases.

(1) Let $X=M(j),\,j>0$. By Theorem 5(4) and Lemmas 8 and 9, $C(X)\cap I(j-1)=L(j-1)$ and $C(X)\cap I(k)=0$ for $k\not= j-1,j$. On the other hand whenever $D\cap I(k)\not= 0$,  $D\supset\text{soc}\,I(k)=L(k)$. But $L(k)$ is in the coradical of $D$ iff $k=j-1,j$. Thus we reduce to two cases

(i) Suppose $D\cap I(j)\supsetneqq M(j)$, then $D\supseteqq W(j+1)$. Since $L(j+1)$ is a composition factor of $W(j+1)$, but doesn't lie in the $\text{corad}(D)$, we arrive at a contradiction.

(ii) Assume $D\cap I(j-1)\supsetneqq L(j-1)$. Then by the same argument as in (i), $D\cap I(j-1)\subseteqq W(j)$. Suppose we have equality there. Using Lemma 11 we would have $NW(j)=0$. Now $N=\frak{m_j m_{j-1}}$ and $\frak{m}_{j-1}L(j)=L(j)$ on account of $\frak{m}_{j-1}+\frak{m}_j=U$. But $L(j)$ is the top composition factor of $W(j)$, hence $NW(j)=\frak{m}_jW(j)=0$, which is impossible, since $W(j)$ is not semisimple. This completes case (1).

(2) Let $X=I(j)$ for a $j\not= 1$. Combining Theorem 5 and Lemmas 8 and 9 we get $C(X)\cap I(j-1)=W(j)$ (for $j>0$) and $C(X)\cap I(j+1)=M(j+1)$.

On the other hand $L(k)$ is in the socle of $D$ iff $k\in \{j-1,j,j+1\}$. Therefore $D\cap I(k)\not= 0$ for those $k$ only. Were $D\cap I(k)\supsetneqq C(X)\cap I(k)$ then $D$ would have a simple subcomodule $L(s)$ with $s\notin\{ j-1,j,j+1\}$, a contradiction.

The last case to consider is when $X=I(1)$. Here we must dispose of the possibility $D\cap I(0)=I(0)$. If this holds, then $NI(0)=0$ follows. Now by definition $N=\frak{m}_1\frak{m}_0\frak{m}_2\frak{m}_1$; hence $NI(0)=(\frak{m}_1\frak{m}_0)I(0)$ on account of $\frak{m}_iL(0)=L(0)$ for $i=1,2$. Further, $\frak{m}_0I(0)=W(1)$, for otherwise $I(0)/L(0)$ would be a semisimple $U$- module, a contradiction. It follows that $\frak{m}_1W(1)=0$, which is again impossible, since $W(1)$ is not semisimple. We conclude that $D\cap I(0)=W(1)$, and the proof is complete
\qed
\enddemo

\remark{Remark \rom{4}} The annihilator in case (1) in the above proof is a particular instance of the annihilator of an $R$- module $E$ in $\text{Ext}_R^1(Y,X)$ where $X$ and $Y$ are simple finite- dimensional $R$- modules. Denote by $\frak{m}_X$ and $\frak{m}_Y$ the annihilators in $R$ of those simple modules. The following holds: $\text{ann}_RE=\frak{m}_X\frak{m}_Y $ if and only if $\dim\text{Ext}_R^1(Y,X)\le 1$. This is the case, of course, for every two $U$- simples.
\endremark

\enddocument